\documentclass[10pt]{amsart}

\usepackage{amsmath,amsthm,amsfonts,amscd,amssymb}
\usepackage{hyperref,wasysym}
\usepackage{verbatim}

\newtheorem{thm}{Theorem}[section]
\newtheorem{cor}[thm]{Corollary}
\newtheorem{lem}[thm]{Lemma}

\theoremstyle{definition}

\theoremstyle{remark}
\newtheorem{rem}{Remark}[section]

\begin{document}

\title{On zeros of multilinear polynomials}

\author{Maxwell Forst}
\author{Lenny Fukshansky}\thanks{Fukshansky was partially supported by the Simons Foundation grant \#519058}

\address{Institute of Mathematical Sciences, Claremont Graduate University, Claremont, CA 91711}
\email{maxwell.forst@cgu.edu}
\address{Department of Mathematics, 850 Columbia Avenue, Claremont McKenna College, Claremont, CA 91711}
\email{lenny@cmc.edu}


\begin{abstract} We consider multivariable polynomials over a fixed number field, linear in some of the variables. For a system of such polynomials satisfying certain technical conditions we prove the existence of search bounds for simultaneous zeros with respect to height. For a single such polynomial, we prove the existence of search bounds with respect to height for zeros lying outside of a prescribed algebraic set. We also obtain search bounds in the case of homogeneous multilinear polynomials, which are related to a so-called ``sparse" version of Siegel's lemma. Among the tools we develop are height inequalities that are of some independent interest.
\end{abstract}

\maketitle

\def\A{{\mathcal A}}
\def\B{{\mathcal B}}
\def\C{{\mathcal C}}
\def\D{{\mathcal D}}
\def\F{{\mathcal F}}
\def\x{{\mathcal H}}
\def\I{{\mathcal I}}
\def\J{{\mathcal J}}
\def\K{{\mathcal K}}
\def\L{{\mathcal L}}
\def\M{{\mathcal M}}
\def\N{{\mathcal N}}
\def\O{{\mathcal O}}
\def\R{{\mathcal R}}
\def\s{{\mathcal S}}
\def\V{{\mathcal V}}
\def\W{{\mathcal W}}
\def\X{{\mathcal X}}
\def\Y{{\mathcal Y}}
\def\H{{\mathcal H}}
\def\Z{{\mathcal Z}}
\def\OO{{\mathcal O}}
\def\BB{{\mathbb B}}
\def\cee{{\mathbb C}}
\def\pee{{\mathbb P}}
\def\que{{\mathbb Q}}
\def\real{{\mathbb R}}
\def\zed{{\mathbb Z}}
\def\hyp{{\mathbb H}}
\def\aa{{\mathfrak a}}
\def\HH{{\mathfrak H}}
\def\qbar{{\overline{\mathbb Q}}}
\def\eps{{\varepsilon}}
\def\ahat{{\hat \alpha}}
\def\bhat{{\hat \beta}}
\def\gt{{\tilde \gamma}}
\def\h{{\tfrac12}}
\def\be{{\boldsymbol e}}
\def\bei{{\boldsymbol e_i}}
\def\bff{{\boldsymbol f}}
\def\ba{{\boldsymbol a}}
\def\bb{{\boldsymbol b}}
\def\bc{{\boldsymbol c}}
\def\bm{{\boldsymbol m}}
\def\bk{{\boldsymbol k}}
\def\bi{{\boldsymbol i}}
\def\bl{{\boldsymbol l}}
\def\bq{{\boldsymbol q}}
\def\bu{{\boldsymbol u}}
\def\bt{{\boldsymbol t}}
\def\bs{{\boldsymbol s}}
\def\bv{{\boldsymbol v}}
\def\bw{{\boldsymbol w}}
\def\bx{{\boldsymbol x}}
\def\bX{{\boldsymbol X}}
\def\bz{{\boldsymbol z}}
\def\bwy{{\boldsymbol y}}
\def\bY{{\boldsymbol Y}}
\def\bL{{\boldsymbol L}}
\def\baa{{\boldsymbol\alpha}}
\def\bbb{{\boldsymbol\beta}}
\def\bet{{\boldsymbol\eta}}
\def\bxi{{\boldsymbol\xi}}
\def\bo{{\boldsymbol 0}}
\def\bol{{\boldkey 1}_L}
\def\ep{\varepsilon}
\def\p{\boldsymbol\varphi}
\def\q{\boldsymbol\psi}
\def\rank{\operatorname{rank}}
\def\aut{\operatorname{Aut}}
\def\lcm{\operatorname{lcm}}
\def\sgn{\operatorname{sgn}}
\def\spn{\operatorname{span}}
\def\md{\operatorname{mod}}
\def\Norm{\operatorname{Norm}}
\def\dim{\operatorname{dim}}
\def\det{\operatorname{det}}
\def\Vol{\operatorname{Vol}}
\def\rk{\operatorname{rk}}
\def\Gal{\operatorname{Gal}}
\def\WR{\operatorname{WR}}
\def\WO{\operatorname{WO}}
\def\GL{\operatorname{GL}}
\def\pr{\operatorname{pr}}

\section{Introduction and statement of results}
\label{intro}

Hilbert's 10th problem asks for an algorithm to decide whether a given Diophantine equation has an integer solution. By a celebrated result of Matiyasevich~\cite{matiyasevich}, such an algorithm does not exist in general. On the other hand, for linear Diophantine equations solutions are classically given by the Euclidean algorithm. Further, there are also known algorithms for quadratic polynomials (see, for instance~\cite{segal}). An important approach to the problem of finding such algorithms for different classes of polynomials is through the use of {\it search bounds}, as described in~\cite{masser1}. Suppose we can prove that a given equation has integer solutions if and only if it has a solution of norm bounded by some explicit function of the coefficients of this equation. Then a search through a finite set of all integer points with norm bounded by this function provides an algorithm that decides whether a solution exists and finds at least one such solution if it exists. 

In fact, search bounds for zeros of polynomial equations have been studied quite extensively over more general rings and fields as well: in these more general situation the role of norm guaranteeing the finiteness of a searchable set is played by an appropriate {\it height function} (we recall the definitions and properties of related height functions $H$, $\H$ and $h$ in Section~\ref{notation}). Search bounds with respect to height over different fields and rings of arithmetic interest for a collection of linear forms are given by various forms of Siegel's lemma: we recall it in Section~\ref{siegel} and prove a certain ``sparse" version of it (Theorem~\ref{siegel_sparse}). The subject of search bounds for quadratic polynomials has been started by a classical theorem of Cassels~\cite{cassels_qf}; see~\cite{lf:cassels} for a detailed overview of a large body of work on various extensions and generalizations of this important theorem. Additionally, there are search bounds for integral cubic forms in a sufficient number of variables~\cite{cubic}, as well as for systems of integral forms under certain technical non-singularity conditions~\cite{qr_dioph}.

The main goal of this paper is to investigate bounds on height of ``small" solutions to polynomial equations, linear in some of the variables, over number fields. From here on, let $K$ be a number field of degree $d$ and discriminant $\Delta_K$. Let $F(x_1,\dots,x_n) \in K[x_1,\dots,x_n]$ and let $1 \leq k < n$. Let $[n] = \{1,\dots,n\}$, $I = \{i_1,\dots,i_k\} \subset [n]$, and $I' = [n] \setminus I$. Let $\bx_{I'} = (x_j)_{j \in I'}$. We will say that $F$ is linear in $I$-separated variables if
\begin{equation}
\label{I_sep}
F(x_1,\dots,x_n) = \sum_{j=1}^k x_{i_j} F_j(\bx_{I'}) + F_{k+1}(\bx_{I'}),
\end{equation}
where $F_j(\bx_{I'}) \in K[\bx_{I'}]$ for $1 \leq j \leq k+1$ are any polynomials in $n-k$ variables indexed by $I'$ with coefficients in $K$. For a polynomial $F(x_1,\dots,x_n) \in K[x_1,\dots,x_n]$, we define its zero-set over~$K$
$$Z_K(F) = \left\{ \bz \in K^n : F(\bz) = 0 \right\}.$$
We also write $\N(F)$ for the number of nonzero monomials of $F$. With this notation, we can state our first main result.

\begin{thm} \label{multi_zeros} Let $I$ be as above and let 
$$F_l(x_1,\dots,x_n) = \sum_{j=1}^k x_{i_j} F_{l,j}(\bx_{I'}) + F_{l,k+1}(\bx_{I'}),\ 1 \leq l \leq k$$
be polynomials over $K$ of respective degrees $m_1,\dots,m_k$ linear in $I$-separated variables as in~\eqref{I_sep}. Consider the inhomogeneous system
\begin{align}
\label{system-F} 
\left. \begin{array}{ll}
F_1(x_1,\dots,x_n) = \sum_{j=1}^k x_{i_j} F_{1,j}(\bx_{I'}) + F_{1,k+1}(\bx_{I'}) &= 0 \\
& \vdots \\
F_k(x_1,\dots,x_n) = \sum_{j=1}^k x_{i_j} F_{k,j}(\bx_{I'}) + F_{k,k+1}(\bx_{I'}) &= 0 
\end{array}
\right\}
\end{align}
of linear equations in the variables $x_{i_1},\dots,x_{i_k}$ with coefficients $F_{l,j}(\bx_{I'})$, $1 \leq l \leq k$, $1 \leq j \leq k+1$. Assume that the matrix $\F := (F_{l,j}(\bx_{I'}))_{1 \leq l \leq k, 1 \leq j \leq k}$ of the corresponding homogeneous system has the same rank as the coefficient matrix of inhomogeneous system, i.e., $\F$ augmented by the column $(F_{l,k+1}(\bx_{I'}))_{1 \leq l \leq k}$. Then $\bigcap_{l=1}^k Z_K(F_l) \neq \emptyset$ and there exists a point $\bz \in \bigcap_{l=1}^k Z_K(F_l)$ with
$$h(\bz) \leq k^{k+1} |\Delta_K|^{\frac{1}{d}} \left( \frac{D+2}{2} \right)^{2km+1} \left( \N \HH \right)^{2k},$$
where 
$$D = \sum_{l=1}^k m_l,\ m = \max_{1 \leq l \leq k} m_l,$$
$$\N = \max_{1 \leq l \leq k} \N(F_l),\ \HH = \max_{1 \leq l \leq k} h(F_l).$$
\end{thm}

\noindent
We prove this theorem in Section~\ref{many_poly}, where we also show how our method of proof leads to an explicit algorithm for finding a simultaneous zero of the polynomial system in question. Our main tools are the Bombieri--Vaaler Siegel's lemma (Theorem~\ref{siegel_lem}), a non-vanishing lemma for polynomials related to Alon's Combinatorial Nullstellensatz (Lemma~\ref{nonvanish}) and a collection of height inequalities that we discuss in Section~\ref{notation}. One of these height inequalities that we prove, a bound on the height of the inverse of a nonsingular matrix (Lemma~\ref{inverse}), is of some independent interest and may have other applications in Diophantine geometry.  
\smallskip

In the case of a single polynomial, we can prove a similar result but with additional avoidance conditions.

\begin{thm} \label{FV_null} Let $n \geq 2$ and $F(\bx)$ and $P(\bx)$ be polynomials in $n$ variables over $K$ of degrees $g$ and $m$, respectively. Assume further that $F$ is linear in at least one of the variables and $Z_K(F) \not\subseteq Z_K(P)$. Then there exists a point $\bz \in Z_K(F) \setminus Z_K(P)$ such that
\begin{equation}
\label{F_null_bnd-1}
h(\bz) \leq \N(F) \left( \frac{m (2g-1) + 2}{2} \right)^{g+1} h(F),
\end{equation}
where $\N(F)$ is the number of monomials of $F$.
\end{thm}

\noindent
We prove this theorem in Section~\ref{V_poly}. Finally, we separately discuss the case of homogeneous multilinear polynomials. We refer to a homogeneous polynomial $F(x_1,\dots,x_n) \in K[x_1,\dots,x_n]$ of degree $g$ but linear in every variable as an $(n,g)$-multilinear form over~$K$. Such forms have many zeros; in particular, they vanish on all sufficiently sparse vectors, specifically on vectors with no more than $g-1$ nonzero coordinates. We use this observation to obtain the following search bounds which we prove in Section~\ref{ml_forms}.

\begin{thm} \label{mf_zeros}  Let $V \subseteq K^n$ be an $m$-dimensional subspace and $F$ a multilinear $(n,g)$-form over $K$. Assume that $m+g-1 > n$ and $g > 1$. Then $V$ contains a basis $\bx_1,\dots,\bx_m$ of vectors such that $F(\bx_1) = \dots = F(\bx_m) = 0$ and
$$H(\bx_i) \leq \sqrt{2} m |\Delta_K|^{\frac{m+1}{2d}} \H(V).$$
for each $1 \leq i \leq m$. Further, suppose that $P(x_1,\dots,x_n) \in K[x_1,\dots,x_n]$ is a polynomial such that the dimension of the subspace of $V$ that $Z_K(P,V) := Z_K(P) \cap V$ spans is
$$D(P,V) := \dim \spn_K Z_K(P,V) < m.$$
Then there exists a point $\bz \in V \setminus Z_K(P,V)$ such that $F(\bz) = 0$ and
$$H(\bz) \leq \sqrt{2} m |\Delta_K|^{\frac{m+1}{2d}} \H(V).$$
\end{thm}
\smallskip

\noindent
This theorem is a direct consequence of our ``sparse" Siegel's lemma (Theorem~\ref{siegel_sparse}). For completeness purposes, we present the bounds analogous to Theorem~\ref{mf_zeros} in the simple case of linear forms ($g=1$) separately in Lemma~\ref{linear}. We are now ready to proceed.
\bigskip

\section{Notation and heights}
\label{notation}

Let $K$ be a number field of degree $d := [K : \que] \geq 1$ and let $r_1, r_2$ be the numbers of real and conjugate pairs of complex embeddings, respectively, so that
$$d = r_1 + 2r_2.$$
Let $\Delta_K$ be the discriminant of $K$, and write $M(K)$ for the set of places of $K$. For each $v \in M(K)$ let $d_v = [K_v : \que_v]$ be the local degree, then for each $u \in M(\que)$, $\sum_{v \mid u} d_v = d$. We normalize the absolute values so that for each nonzero $x \in K$ the product formula reads
$$\prod_{v \in M(K)} |x|_v^{d_v} = 1.$$
Let $n \geq 2$, and for any place $v \in M(K)$ and $\bx = (x_1,\dots,x_n) \in K^n$ define the corresponding sup-norm
$$|\bx|_v = \max \{ |x_1|_v,\dots,|x_n|_v \}.$$
If $v \mid \infty$, we also define Euclidean norm 
$$\|\bx\|_v = \left( \sum_{i=1}^n |x_i|_v^2 \right)^{1/2}.$$
Then we define two projective height functions $H, \H : K^n \to \real_{\geq 0}$ as follows:
$$H(\bx) = \left( \prod_{v \in M(K)} |\bx|_v^{d_v} \right)^{1/d},\ \H(\bx) = \left( \prod_{v \nmid \infty} |\bx|_v^{d_v} \times \prod_{v \mid \infty} \|\bx\|_v^{d_v} \right)^{1/d}.$$
These heights are absolute, meaning that they are the same when computed over any number field $K$ containing coordinates of $\bx$: this is due to the normalizing exponent~$1/d$ in the definition. Further, for each nonzero $\bx \in K^n$,
\begin{equation}
\label{ht_ineq_sqrt}
1 \leq H(\bx) \leq \H(\bx) \leq \sqrt{n} H(\bx).
\end{equation}
The inhomogeneous (Weil) height $h : K^n \to \real_{\geq 1}$ is defined by
$$h(\bx) = H(1,\bx) \geq H(\bx)$$
for every $n \geq 1$, thus including a height on algebraic numbers. We also write $H(A)$ and $h(A)$ for the projective and inhomogeneous heights, respectively, of an $m \times n$ matrix $A$ viewed as a vector in $K^{mn}$ and $H(F)$, $h(F)$ for the respective heights of the coefficient vector of a polynomial $F$ over $K$.

We can also define Schmidt height on the subspaces of $K^n$ as follows. Let $V = \spn_K \{ \bx_1,\dots,\bx_m \}$ be an $m$-dimensional subspace of $K^n$, $1 \leq m \leq n$, with a basis $\bx_1,\dots,\bx_m$. The wedge product of these basis vectors $\bx_1 \wedge \dots \wedge \bx_m$ can be viewed as a vector in $K^{\binom{n}{m}}$ under lexicographic embedding: this is a vector of Grassmann (or Pl\"ucker) coordinates of $V$, which is uniquely defined up to a constant multiple. Define
$$\H(V) = \H(\bx_1 \wedge \dots \wedge \bx_m).$$
The product formula guarantees that this definition is independent of the choice of a basis for $V$. We also define Schmidt height on matrices: for an $n \times m$ matrix $A$ over $K$, $1 \leq m \leq n$ we let
$$\H(A) = \H\left( \spn_K \left\{ \ba_1,\dots,\ba_m \right\} \right),$$
where $\ba_1,\dots,\ba_m$ are column vectors of $A$. If $m > n$, we define $\H(A)$ to be $\H(A^{\top})$. Suppose the $m$-dimensional vector subspace $V \subset K^n$ is described as
$$V = \left\{ A \bx : \bx \in K^m \right\} = \left\{ \bwy \in K^n : B \bwy = \bo \right\}$$
for the $n \times m$ matrix $A$ and $(n-m) \times n$ matrix $B$ over $K$, respectively. Then the Brill--Gordan duality principle~\cite{gordan} (also see Theorem~1 on page 294 of~\cite{hodge_pedoe}) states that
\begin{equation}
\label{duality_p}
\H(A) = \H(B) = \H(V).
\end{equation}
 
\smallskip

We also review here several properties of height functions that we will need. The first one is a well known fact, which is an immediate corollary of Theorem~1 of~\cite{struppeck_vaaler}.

\begin{lem} \label{intersect} Let $U_1$ and $U_2$ be subspaces of $K^n$. Then
$$\H(U_1 \cap U_2) \leq \H(U_1) \H(U_2).$$
\end{lem}

\noindent
The second one is Lemma 4.7 of \cite{roy_thunder}.

\begin{lem} \label{lem_4.7} Let $V$ be an $m$-dimensional subspace of $K^n$, $1 \leq m \leq n$, and let $\bx_1,\dots,\bx_m$ be any basis for $V$. Then
$$\H(V) \leq \H(\bx_1) \cdots \H(\bx_m) .$$
\end{lem}

\noindent
The next lemma gives a bound on the height $H$ of the inverse of a square nonsingular matrix (viewed as a vector) in terms of the height of the matrix itself.

\begin{lem} \label{inverse} Let $A \in \GL_n(K)$, then
\begin{equation}
\label{inverse_bnd}
H(A^{-1}) \leq \left( \sqrt{n} H(A) \right)^{n-1},
\end{equation}
and
\begin{equation}
\label{inverse_bnd-1}
h(A^{-1}) \leq n^n |\Delta_K|^{\frac{1}{d}} h(A)^{2n-1}. 
\end{equation}
\end{lem}

\proof
Let us write $\ba_1,\dots,\ba_n$ for the row vectors of $A$, and for each $1 \leq j \leq n$ let $A_j$ be the $(n-1) \times n$ submatrix of $A$ obtained by deleting the $j$-th row $\ba_j$. For each $j$,
$$V_j = \left\{ \bx \in K^n : A_j \bx = \bo \right\}$$
is a $1$-dimensional subspace. Let $\bx_j$ be any nonzero point in $V_j$ and notice that $\ba_j \bx_j \neq 0$: if it was equal to $0$, then $A\bx_j = \bo$, contradicting the assumption that $A$ is nonsingular. Then take $\bb_j = \frac{1}{\ba_j \bx_j} \bx_j$, then  $\ba_i \bb_j = 0$ for every $i \neq j$ and $\ba_j \bb_j = 1$. Take $B$ to be the $n \times n$ matrix whose columns are the vectors $\bb_1,\dots,\bb_n$, then $AB$ is the identity matrix, so $B = A^{-1}$. Notice also that for every $j$,
\begin{equation}
\label{H-inv}
H(\bb_j) \leq \H(\bb_j) = \H(V_j) = \H(A_j) \leq \prod_{i=1, i \neq j}^n \H(\ba_i) \leq n^{\frac{n-1}{2}} \prod_{i=1, i \neq j}^n H(\ba_i),
\end{equation}
by~\eqref{duality_p} combined with Lemma~\ref{lem_4.7} and~\eqref{ht_ineq_sqrt}. Since $H(B) = \max_{1 \leq j \leq n} H(\bb_j)$ and $H(\ba_i) \leq H(A)$ for every $1 \leq i \leq n$, \eqref{inverse_bnd} follows.

To obtain the second inequality, we will use Siegel's lemma (see Theorem~\ref{siegel_lem} below) to choose a specific vector $\bx_j \in V_j$ such that
\begin{equation}
\label{h-inv-1}
h(\bx_j) \leq \left( \left( \frac{2}{\pi} \right)^{2r_2} |\Delta_K| \right)^{1/2d} \H(V_j) \leq |\Delta_K|^{\frac{1}{2d}} \H(V_j).
\end{equation}
Now take $\bb_j = \frac{1}{\ba_j \bx_j} \bx_j$, then 
\begin{eqnarray}
\label{h-inv-2}
h(\bb_j) & \leq & h \left( \sum_{l=1}^n a_{jl} x_{jl}, \bx_j \right) \leq h \left( \sum_{l=1}^n a_{jl} x_{jl} \right) h(\bx_j) \leq n h(\ba_j) h(\bx_j)^2 \nonumber \\
& \leq & n |\Delta_K|^{\frac{1}{d}} h(\ba_j) \H(V_j)^2,
\end{eqnarray}
where the last inequality follows by~\eqref{h-inv-1}. Combining~\eqref{h-inv-2} with~\eqref{H-inv} and observing that $h(B) = \max_{1 \leq j \leq n} h(\bb_j)$ and $h(\ba_i) \leq h(A)$ for every $1 \leq i \leq n$, we obtain~\eqref{inverse_bnd-1}. 
\endproof

\noindent
Next we present a useful bound on the height of a vector whose coordinates are images of a given point under a collection of polynomials.

\begin{lem} \label{poly_ht} Let $F_1,\dots,F_k$ be polynomials of respective degrees $m_1,\dots,m_k$ in $K[x_1,\dots,x_n]$ and $\bz \in K^n$. Then
\begin{equation}
\label{poly_ht_bnd}
h(F_1(\bz),\dots,F_k(\bz)) \leq \N \HH h(\bz)^m,
\end{equation}
where $\N = \max_{1 \leq i \leq k} \N(F_i)$, $\HH = \max_{1 \leq i \leq k} h(F_i)$ and $m = \max_{1 \leq i \leq k} m_i$.
\end{lem}

\proof
For each $1 \leq i \leq k$, let us write
$$F_i(\bx) = \sum_J f_J \bx^J,$$
where $J = (j_1,\dots,j_n)$ is a multi-index with $0 \leq j_l \leq m_i$ for each $1 \leq l \leq n$. Let $\bz \in K^n$, then for every $v \nmid \infty$ in $M(K)$,
$$|F_i(\bz)|_v \leq \max_J \left\{ |f_J|_v \prod_{l=1}^n |z_l|_v^{j_l} \right\} \leq \max_J \{ |f_J|_v \} \max \{1, |z_1|_v,\dots,|z_n|_v \}^{m_i},$$
and for $v \mid \infty$ in $M(K)$,
$$|F_i(\bz)|_v \leq \N(F_i) \max_J \left\{ |f_J|_v \prod_{l=1}^n |z_l|_v^{j_l} \right\} \leq \N(F_i) \max_J \{ |f_J|_v \} \max \{1, |z_1|_v,\dots,|z_n|_v \}^{m_i}.$$
Then
\begin{eqnarray*}
h(F_1(\bz),\dots,F_k(\bz)) & \leq & \prod_{v \in M(K)} \max \left\{ 1,|F_1(\bz)|_v,\dots,|F_k(\bz)|_v \right\}^{\frac{d_v}{d}} \\
& \leq & \left( \max_{1 \leq i \leq k} \N(F_i) \right) \left( \max_{1 \leq i \leq k} h(F_i) \right) h(\bz)^{\max_{1 \leq i \leq k} m_i}.
\end{eqnarray*}
This is precisely~\eqref{poly_ht_bnd}.
\endproof

\noindent
The next lemma bounds the height of a polynomial under a linear transformation. 

\begin{lem} \label{ht_poly_lin} Let $F(x_1,\dots,x_n) \in K[x_1,\dots,x_n]$ be a polynomial of degree $m$ and let $A$ be an $n \times k$ matrix over $K$. For a variable vector $\bwy = (y_1,\dots,y_k)$, define $G(\bwy) = F(A\bwy^{\top})$, then $G$ is a polynomial in $k$ variables over $K$ of degree $\leq m$. Further,
$$h(G) \leq k^m \N(F) h(F) h(A)^m.$$
\end{lem}

\proof
Write $A = (a_{ij})_{1 \leq i \leq n, 1 \leq j \leq k}$ and notice that
$$A\bwy^{\top} = \left( \sum_{j=1}^k a_{1j} y_j,\dots,\sum_{j=1}^k a_{nj} y_j \right)^{\top}.$$
Then $G(\bwy) = F \left( \sum_{j=1}^k a_{1j} y_j,\dots,\sum_{j=1}^k a_{nj} y_j \right)$. Each linear form $\sum_{j=1}^k a_{ij} y_j$ has coefficient vector of inhomogeneous height $\leq h(A)$, and $F$ has $\N(F)$ monomials each of degree no bigger than $m$ and height no bigger than $h(F)$. Therefore
$$h(G) \leq \N(F) h(F) (kh(A))^m,$$
where $k^m$ is an upper bound on the binomial coefficients which occur when taking $m$-th power of a linear form. This is precisely the result of the lemma.
\endproof

Next, we introduce the notion of ``sparsity" of vectors. We say that a vector $\bx \in K^n$ is {\it $s$-sparse} for some $1 \leq s \leq n$ if $\bx$ has no more than $s$ nonzero coordinates. Let $\be_1,\dots,\be_n$ be the standard basis for $K^n$, and for an indexing subset $I \subseteq \{1,\dots,n\}$ of cardinality $t$ let
\begin{equation}
\label{coord_hyp}
\hyp_I = \spn_K \{ \be_i : i \in I \}
\end{equation}
be the corresponding $t$-dimensional coordinate subspace of $K^n$. Then every vector in $\hyp_I$ is $t$-sparse and $\H(\hyp_I) = 1$.

The following lemma that will be quite important to us is a rigorous form of the basic principle that a polynomial which is not identically zero cannot vanish ``too much". Somewhat different formulations of this principle can be found in~\cite{cassels} (Lemma~1 on p.~261) as well as in the context of N. Alon's celebrated Combinatorial Nullstellensatz~\cite{alon}. The following formulation, which is most convenient for our purposes follows easily from Lemma~2.2 of~\cite{lf:int}.

\begin{lem} \label{nonvanish} Suppose $P(\bx) \in K[x_1,\dots,x_n]$ is a polynomial of degree $m$ which is not identically $0$. There exists a point~$\bz \in \zed^n$ such that $P(\bz) \neq 0$ and 
$$h(\bz) \leq \frac{m+2}{2}.$$
\end{lem}

\bigskip

\section{On Siegel's lemma}
\label{siegel}

Let us recall the classical Siegel's lemma over a number field $K$ as proved in~\cite{bombieri_vaaler} by Bombieri and Vaaler (see also Theorem~1.1 of~\cite{lf:null} for a generalized formulation over different choices of the ground field).

\begin{thm} \label{siegel_lem} Let $V \subseteq K^n$ be an $m$-dimensional subspace, $1 \leq m \leq n$. There exists a basis $\bx_1,\dots,\bx_m$ for $V$ such that
\begin{equation}
\label{siegel1}
\prod_{i=1}^m h(\bx_i) \leq \left( \left( \frac{2}{\pi} \right)^{2r_2} |\Delta_K| \right)^{m/2d} \H(V).
\end{equation}
\end{thm}

Our main goal in this section to establish a ``sparse" version of Siegel's lemma, i.e., existence of a small-height basis of sparse vectors for $V$. One can ask if a sparse basis necessarily exists. Suppose $\bx_1,\dots,\bx_m$ is some basis for $V$ and write $X = (\bx_1\ \dots\ \bx_m)$ for the corresponding $n \times m$ basis matrix, so that $V = XK^m$. Let $Y$ be the reduced column-echelon form of $X$. Then $V = YK^m$ so column vectors of $Y$ again form a basis for $V$. Further, the first column of $Y$ has the fewest number of zero coordinates, and this number is at least equal to the number of remaining columns, which is $m-1$, since all the entries to the left of the leading $1$ in each column are $0$'s. Hence the column vectors of $Y$ provide a basis of $(n-m+1)$-sparse vectors for $V$. On the other hand, it is not so easy to produce a bound on the height of these basis vectors, since it is tricky to keep track of the coefficients used in the column-echelon reduction. To do this, we need a certain version of Siegel's lemma with avoidance conditions established in~\cite{lf:null} (Corollary~1.5).

\begin{thm} \label{siegel_avoid} Let $V \subseteq K^n$ be an $m$-dimensional subspace, $1 \leq m \leq n$, and $U_1,\dots,U_k \subset K^n$ be subspaces such that $V \not\subseteq \bigcup_{i=1}^k U_i$. Then there exists $\bx \in V \setminus \left( \bigcup_{i=1}^k U_i \right)$ such that
$$h(\bx) \leq \sqrt{2} m |\Delta_K|^{\frac{m+1}{2d}} k^{\frac{1}{d}} \H(V).$$
\end{thm}

With this in mind, we can prove the following result.

\begin{thm} \label{siegel_sparse} Let $V \subseteq K^n$ be an $m$-dimensional subspace, $1 \leq m \leq n$. Then $V$ contains a basis $\bx_1,\dots,\bx_m$ of $(n-m+1)$-sparse vectors such that
\begin{equation}
\label{ss_bnd}
h(\bx_i) \leq \sqrt{2} m |\Delta_K|^{\frac{m+1}{2d}} \H(V).
\end{equation}
for each $1 \leq i \leq m$.
\end{thm}

\proof
Take $t=n-m+1$. First suppose that $V$ is contained in some $t$-dimensional coordinate subspace. Then the basis guaranteed by Theorem~\ref{siegel_lem} is $(n-m+1)$-sparse and satisfies the bound~\eqref{ss_bnd}. Hence we can assume that $V$ is not contained in any such $\hyp_I$. Let
$$\J(n,t) = \left\{ I \subset \{1,\dots,n\} : |I| = t \right\}.$$
Then cardinality of $\J(n,t)$ is equal to $\binom{n}{t}$. Pick some $I_1 \in \J(n,t)$ and let $m_1$ be the dimension of the subspace $V \cap \hyp_{I_1} \subsetneq V$. Since $\dim V + \dim \hyp_{I_1} > n$, we must have $1 \leq m_1 < m$. Let $\bx_1$ be a vector of smallest height from a basis for $V \cap \hyp_{I_1}$ guaranteed by Theorem~\ref{siegel_lem}, then $\bx_1$ is an $(n-m+1)$-sparse vector satisfying
\begin{equation}
\label{x1}
h(\bx_1) \leq \left( \left( \frac{2}{\pi} \right)^{2r_2} |\Delta_K| \right)^{\frac{1}{2d}} \H(V \cap \hyp_{I_1})^{\frac{1}{m_1}} \leq |\Delta_K|^{\frac{1}{2d}} \H(V),
\end{equation}
since $\H(V \cap \hyp_{I_1}) \leq \H(V) \H(\hyp_{I_1}) = \H(V)$ by Lemma~\ref{intersect}. Let $V_1 = K \bx_1$. Now, by our column-echelon form argument above, we know that $V$ has a basis of $(n-m+1)$-sparse vectors, hence there must exist some $I_2 \in \J(n,t)$ such that $V \cap \hyp_{I_2} \not\subseteq V_1$. Every vector in $V \cap \hyp_{I_2}$ is again $(n-m+1)$-sparse, so let $\bx_2 \in (V \cap \hyp_{I_2}) \setminus V_1$ be a vector guaranteed by Theorem~\ref{siegel_avoid}, then
\begin{equation}
\label{x2}
h(\bx_2) \leq \sqrt{2} m |\Delta_K|^{\frac{m+1}{2d}} \H(V),
\end{equation}
since $k = 1$ in this case, dimension of $V \cap \hyp_{I_2}$ is $\leq m$ and again $\H(V \cap \hyp_{I_2}) \leq \H(V) \H(\hyp_{I_2}) = \H(V)$. Let $V_2 = \spn_K \{ V_1, \bx_2 \}$, then either $V_2 = V$ or, again, there must exist some $I_3 \in \J(n,t)$ such that $V \cap \hyp_{I_3} \not\subseteq V_2$, and hence we can use Theorem~\ref{siegel_avoid} to obtain $\bx_3 \in (V \cap \hyp_{I_3}) \setminus V_2$ with height bounded the same way as in~\eqref{x2}. Continuing in the same manner $m$ times, we obtain an $(n-m+1)$-sparse basis satisfying~\eqref{ss_bnd}, since the bound of~\eqref{x1} is less than the bound of~\eqref{x2}.
\endproof

\bigskip

\section{Zeros of multiple polynomials}
\label{many_poly}

The main goal of this section is to prove Theorem~\ref{multi_zeros}. Let the notation be as in the statement of the theorem.

\proof[Proof of Theorem~\ref{multi_zeros}]
Let $r$ be the rank of the linear system~\eqref{system-F}, which is the same as the rank of the homogeneous system as in the statement of the theorem. If $r < k$, then $k-r$ equations of~\eqref{system-F} are dependent on the rest of them, and so every solution to the rest of the equations is automatically a solution to the whole system. Since $r$ equations are linearly independent, there must exist some $r \times r$ submatrix of the coefficient matrix of these equations with entries $F_{ij}(\bx_{I'})$ which is nonsingular for some $\bx_{I'} \in K^{n-k}$, and hence the determinant of this matrix is not identically zero as a polynomial in the variables of $\bx_{I'}$. We can set all the $k-r$ variables $x_{i_j}$ not corresponding to the columns of this submatrix equal to $0$, and hence reduce to the case of $r$ equations in $r$ variables. Hence we can assume without loss of generality that $r = k$.

Let us rewrite~\eqref{system-F} as
\begin{equation}
\label{matrix_system}
\F(\bx_{I'}) \bx_I = \bff(\bx_{I'}),
\end{equation}
where $\F(\bx_{I'})$ is the $k \times k$ matrix with entries $F_{ij}(\bx_{I'})$, $\bff(\bx_{I'})$ is the $k$-dimensional column vector with coordinates $-F_{i (k+1)}(\bx_{I'})$ and $\bx_I = (x_{i_1},\dots,x_{i_k})^{\top}$ is the variable vector. Then $\F(\bx_{I'})$ is nonsingular for some choice of $\bx_{I'} \in K^{n-k}$, hence $P(\bx_{I'}) := \det(\F(\bx_{I'}))$ is not identically zero as a polynomial in the variables $\bx_{I'}$. Notice that $\deg(P) \leq \sum_{i=1}^k \deg(F_i) = D$, and hence by Lemma~\ref{nonvanish} there exists a point~$\bz_{I'} \in \zed^{n-k}$ such that $P(\bz_{I'}) \neq 0$ and 
\begin{equation}
\label{ht_zI}
h(\bz_{I'}) \leq \frac{D+2}{2}.
\end{equation}
Plugging in $\bz_{I'}$ for $\bx_{I'}$ into~\eqref{matrix_system}, we obtain a nonsingular linear system of $k$ equations in $k$ variables, and hence have a unique solution $\bz_I = \F(\bz_{I'})^{-1} \bff(\bz_{I'})$. Combining $\bz_I$ with $\bz_{I'}$ into appropriately indexed coordinates, we obtain a vector $\bz \in \bigcap_{j=1}^k Z_K(F_j)$ and $H(\bz) = H(\bz_I, \bz_{I'})$. 

We now estimate the height of~$\bz_I$. First notice that, by Lemma~\ref{inverse},
\begin{equation}
\label{htF1}
h(\F(\bz_{I'})^{-1}) \leq k^k |\Delta_K|^{\frac{1}{d}} h(\F(\bz_{I'}))^{2k-1}.
\end{equation}
On the other hand,
\begin{equation}
\label{htF2}
h(\F(\bz_{I'})) = H(1,F_{11}(\bz_{I'}),\dots,F_{kk}(\bz_{I'})) \leq \N h(\bz_{I'})^m \max_{1 \leq i \leq k} h(F_i),
\end{equation}
as well as 
\begin{equation}
\label{htF3}
h(\bff (\bz_{I'})) = H(1,-F_{1(k+1)}(\bz_{I'}),\dots,-F_{k(k+1)}(\bz_{I'})) \leq \N h(\bz_{I'})^m \max_{1 \leq i \leq k} h(F_i),
\end{equation}
both by Lemma~\ref{poly_ht}. Combining~\eqref{htF1}, \eqref{htF2}, \eqref{htF3} and~\eqref{ht_zI}, we obtain:
\begin{eqnarray}
\label{htF4}
h(\bz_I) & = & h \left( \F(\bz_{I'})^{-1} \bff(\bz_{I'}) \right) \leq k h \left(\F(\bz_{I'})^{-1}) h(\bff (\bz_{I'})) \right) \nonumber \\
& \leq & k^{k+1} |\Delta_K|^{\frac{1}{d}} h(\F(\bz_{I'}))^{2k-1} h(\bff (\bz_{I'})) \nonumber \\
& \leq & k^{k+1} |\Delta_K|^{\frac{1}{d}} \left( \N h(\bz_{I'})^m \max_{1 \leq i \leq k} h(F_i) \right)^{2k} \nonumber \\
& \leq & k^{k+1} |\Delta_K|^{\frac{1}{d}} \left( \N \left( \frac{D+2}{2} \right)^m \max_{1 \leq i \leq k} h(F_i) \right)^{2k}.
\end{eqnarray}
Now notice that
\begin{eqnarray*}
h(\bz) & = & H(1,\bz_I, \bz_{I'}) = \prod_{v \in M(K)} \max \left\{ 1, |\bz_I|_v, |\bz_{I'}|_v \right\}^{\frac{d_v}{d}} \\
& \leq & \prod_{v \in M(K)} \left( \max \left\{ 1, |\bz_I|_v \right\}^{\frac{d_v}{d}} \max \left\{ 1, |\bz_{I'}|_v \right\}^{\frac{d_v}{d}} \right) = h(\bz_I) h(\bz_{I'}),
\end{eqnarray*}
and hence the theorem follows from~\eqref{ht_zI} and~\eqref{htF4}.
\endproof

\begin{rem} \label{algorithm} Our argument leads to an algorithm for finding a simultaneous zero $\bz$ of the polynomial system~\eqref{system-F} under the assumption that it exists:

\begin{enumerate}

\item Compute $P(\bx_{I'}) = \det(\F(\bx_{I'}))$ as a polynomial in the variables $\bx_{I'}$.

\item Search through the set of integer points of sup-norm $\leq \deg(P)$ to find $\bz_{I'}$ such that $P(\bz_{I'}) \neq 0$.

\item For this choice of $\bz_{I'}$, compute $\F(\bz_{I'})$ and $\bff(\bz_{I'})$.

\item Compute $\F(\bz_{I'})^{-1}$.

\item Compute $\bz_I = \F(\bz_{I'})^{-1} \bff(\bz_{I'})$.

\item Combine $\bz_I$ and $\bz_{I'}$ according to the indices of the coordinates to obtain the vector $\bz$.

\end{enumerate}
\end{rem}

\bigskip

\section{Zeros of one polynomial}
\label{V_poly}

In this section, we prove existence of a zero of bounded height for a polynomial $F$, linear in at least one variable, outside of an algebraic set not containing the entire zero locus of~$F$. Our main goal is to prove Theorem~\ref{FV_null}.

\proof[Proof of Theorem~\ref{FV_null}]
Assume without loss of generality that $F$ is linear in $x_n$. Define $\bx' := (x_1,\dots,x_{n-1})$, so $\bx = (\bx',x_n)$. We can write
$$F(x_1,\dots,x_n) = x_n F_1(x_1,\dots,x_{n-1}) + F_2(x_1,\dots,x_{n-1}),$$
where $\deg F_1 \leq g-1$, $\deg F_2 \leq g$ and both of them are polynomials in $n-1$ variables $\bx'$ with at least $F_1$ not identically zero. Then we can describe $Z_K(F)$ as the union $Z^1_K(F) \cup Z^2_K(F)$, where
\begin{align}
\label{F_zero}
Z^1_K(F) & = \left\{ \bz \in K^n : F_1(\bz') = F_2(\bz') = 0 \right\},\nonumber \\
Z^2_K(F) & = \left\{ \bz \in K^n : F_1(\bz') \neq0, z_n = -F_2(\bz')/F_1(\bz') \right\}.
\end{align}
Let $Z^2_K(F)' = \{ \bz' : \bz \in Z^2_K(F) \}$, and define
$$Q(\bx') = P \left( x_1,\dots, x_{n-1}, -\frac{F_2(\bx')}{F_1(\bx')} \right) F_1(\bx')^m.$$
Since the degree of $P$ is $m$, $Q$ is a polynomial in $n-1$ variables $\bx'$ on $Z^2_K(F)'$. Let us show that $Q$ is not identically $0$. Suppose it is, then $P(Z^2_K(F)) = 0$. For each $\bz' \in Z^2_K(F)'$ view $P(\bz',x_n)$ as a polynomial in one variable $x_n$. Since 
$$x_n = - F_2(\bz')/F_1(\bz')$$
is a root of this polynomial, the linear factor $x_n+F_2(\bz')/F_1(\bz')$ must divide it, and since we are working over a field, we can say that $F_1(\bz')x_n + F_2(\bz')$ divides $P(\bz',x_n)$ for every $\bz' \in Z^2_K(F)'$. This implies that $P(\bx)$ must be divisible by $x_nF_1(\bx')+F_2(\bx') = F(\bx)$, and therefore $Z_K(F) \subseteq Z_K(P)$. This is a contradiction, and hence there exists $\bz' \in Z^2_K(F)$ such that $Q(\bz') \neq 0$. We want to find such $\bz'$ of bounded height. Notice that 
$$\deg Q \leq m \deg F_2 + m \deg F_1 \leq m(g + g - 1) = m (2g-1).$$
By Lemma~\ref{nonvanish}, there exists $\bz' \in \zed^{n-1}$ such that $Q(\bz') \neq 0$ and
\begin{equation}
\label{Hz}
h(\bz') \leq \frac{\deg Q+1}{2} = \frac{m (2g-1) + 2}{2}.
\end{equation}
Then we can estimate the height of the corresponding point $\bz = \left( \bz', -\frac{F_2(\bz')}{F_1(\bz')} \right)$. Notice that
\begin{eqnarray}
\label{HF_bnd}
h(\bz) & \leq & H(1,F_1(\bz') \bz',F_2(\bz')) \nonumber \\
& = & \prod_{v \in M(K)} \max \left\{ 1,|F_1(\bz') z_1|_v,\dots, |F_1(\bz') z_{n-1}|_v,|F_2(\bz')|_v \right\}^{\frac{d_v}{d}} \nonumber \\
& \leq & \prod_{v \in M(K)} \left( \max \left\{ 1,|F_1(\bz')|_v,|F_2(\bz')|_v \right\} \max\{ 1, |z_1|_v,\dots,|z_{n-1}|_v \} \right)^{\frac{d_v}{d}} \nonumber \\
& \leq & H \left(1,F_1(\bz'),F_2(\bz') \right) h(\bz') \leq \N(F) h(F) h(\bz')^g h(\bz'),
\end{eqnarray}
where the last inequality follows by Lemma~\ref{poly_ht}. Combining~\eqref{Hz} and \eqref{HF_bnd} yields~\eqref{F_null_bnd-1}.
\endproof

\begin{rem} Notice that description~\eqref{F_zero} immediately implies that $Z_K(F)$ is not empty, and therefore taking $P$ to be a nonzero constant polynomial we see that $F$ has a zero $\bz \in K^n$ with $h(\bz) \leq \N(F) h(F)$. Further, our argument allows for an explicit algorithm to find the point~$\bz$ in question, similar to the procedure described in Remark~\ref{algorithm} for Theorem~\ref{multi_zeros}: here, we just need to do a finite search for an integer point $\bz'$ so that $Q(\bz') \neq 0$ and define $\bz = \left( \bz', -\frac{F_2(\bz')}{F_1(\bz')} \right)$.
\end{rem}

One can further ask if a result similar to Theorem~\ref{FV_null} holds with a restriction to a subspace $V$ of $K^n$. The problem here is that the restriction of our polynomial $F$ to $V$ may no longer be linear in any of the variables. For example, if $F(x_1,x_2) = x_1x_2$ and $V = K \begin{pmatrix} 1 \\ 1 \end{pmatrix} \subset K^2$, then the restriction of $F$ to $V$ is
$$F_V(x) = F(x,x) = x^2,$$
and hence is not linear. On the other hand, we can prove a simple lemma in case $\dim_K V = 1$.

\begin{lem} \label{FV_null-1} Let $F$, $P$ be as in Theorem~\ref{FV_null} and suppose $V$ is a one-dimensional subspace of $K^n$ such that
$$Z_K(F,V) := Z_K(F) \cap V \not\subseteq Z_K(P).$$
Then there exists $\bz \in Z_K(F,V) \setminus Z_K(P)$ such that
$$h(\bz) \leq \left( \left( \frac{2}{\pi} \right)^{2r_2} |\Delta_K| \right)^{\frac{m+1}{2d}} \N(F)^{\frac{3}{2}} h(F) \H(V)^{m+1}.$$
\end{lem}

\proof
Now suppose that $V$ is a one-dimensional subspace of $K^n$, i.e., $\ell = 1$. Then $V = K\bwy$ for some vector $\bwy \in K^n$ and $F(\alpha \bwy) = 0$ for some $\alpha \in K$ such that $\bz := \alpha \bwy \not\in Z_K(P)$, since $Z_K(F,V) \not\subseteq Z_K(P)$. By Theorem~\ref{siegel_lem}, we can choose $\bwy$ such that
\begin{equation}
\label{y_siegel}
h(\bwy) \leq  \left( \left( \frac{2}{\pi} \right)^{2r_2} |\Delta_K| \right)^{\frac{1}{2d}} \H(V).
\end{equation}
Then $F_V$ is a polynomial in one variable of degree $\leq m$ with
\begin{equation}
\label{FV5.1}
h(F_V) \leq \N(F) h(F) h(\bwy)^m \leq \left( \left( \frac{2}{\pi} \right)^{2r_2} |\Delta_K| \right)^{\frac{m}{2d}} \N(F) h(F) \H(V)^m,
\end{equation}
by Lemma~\ref{ht_poly_lin} combined with~\eqref{y_siegel}. Let $\alpha_1,\dots,\alpha_k$ be the roots of $F_V$, $k \leq m$, with repetition if necessary. Then Lemma~2 of~\cite{vaaler_pinner} combined with~\eqref{ht_ineq_sqrt} guarantees that
\begin{equation}
\label{mahler}
\prod_{i=1}^k h(\alpha_i) \leq \sqrt{\N(F)} h(F_V).
\end{equation}
Observing that $\max_{1 \leq i \leq k} h(\alpha_i) \leq \prod_{i=1}^k h(\alpha_i)$ and combining~\eqref{FV5.1} with~\eqref{mahler}, we obtain
$$\max_{1 \leq i \leq k} h(\alpha_i \bwy) \leq h(\bwy) \max_{1 \leq i \leq k} h(\alpha_i) \leq \left( \left( \frac{2}{\pi} \right)^{2r_2} |\Delta_K| \right)^{\frac{m+1}{2d}} \N(F)^{\frac{3}{2}} h(F) \H(V)^{m+1},$$
which is the bound of the lemma. This completes the proof.
\endproof
\bigskip

\section{Zeros of multilinear forms}
\label{ml_forms}

In this section we prove Theorem~\ref{mf_zeros}. For $n \geq 2$, let $n \geq g \geq 1$ and let $F(x_1,\dots,x_n) \in K[x_1,\dots,x_n]$ be a multilinear $(n,g)$-form, which is not identically zero.

\begin{lem} \label{Z1} There exists a nonzero point $\bz \in K^n$ such that $F(\bz) = 0$ and
\begin{equation}
\label{ht_Z1}
H(\bz) \leq H(F).
\end{equation}
\end{lem}

\proof
We argue by induction on $n$. If $n=2$, then
$$F(x_1,x_2) = ax_1+bx_2, \text{ or } F(x_1,x_2) = cx_1x_2,$$
for some $a,b,c \in K$. In the first case, $F(-b,a) = 0$, and in the second $F(0,1)=0$. In either case, the nontrivial zero $\bz = (-b,a)$ or $\bz = (0,1)$ satisfies the bound~\eqref{ht_Z1}.

Suppose now $n > 2$. If $n=g$, then
$$F(x_1,\dots,x_n) = cx_1 \cdots x_n,$$
and so $F(0,\dots,0,1)=0$. Hence assume $n > g$, then for some $1 \leq i \leq n$ we can write
$$F(x_1,\dots,x_n) = x_i F_1(\bx'_i) + F_2(\bx'_i),$$
where
$$\bx'_i = (x_1,\dots,x_{i-1},x_{i+1},\dots,x_n)$$
is the vector of $n-1$ variables excluding $x_i$ and $F_1,F_2$ are multilinear forms in $n-1$ variables, not identically zero. By the induction hypothesis, there exists a nonzero point $\bz' \in K^{n-1}$ such that $F_2(\bz') = 0$ satisfying~\eqref{ht_Z1}. Define $\bz$ by inserting $0$ for the $i$-th coordinate in $\bz'$, then $F(\bz) = 0$ and $\bz$ satisfies~\eqref{ht_Z1}.
\endproof

On the other hand, if $V \subseteq K^n$ is a subspace of $K^n$, then $F$ may not necessarily have nontrivial zeros on $V$. Indeed, the form $x_1+x_2$ has no nontrivial zeros on the subspace $\{ (a,2a) : a \in \que \}$ of $\que^2$. There are, however, some situations when $F$ is guaranteed to have nontrivial zeros on $V$, and then we can find such zeroes of small height. 

\begin{lem} \label{sparse_zeros} Assume $g > 1$ and let $F$ be a multilinear $(n,g)$-form over $K$ and $\bx \in K^n$ a $(g-1)$-sparse vector. Then $F(\bx) = 0$.
\end{lem}

\proof
The vector $\bx$ has no more than $g-1$ nonzero coordinates. Hence $F(\bx) = 0$, since every monomial of $F$ has degree $g$ and is linear in each variable, hence is a product of $g > 1$ distinct variables. 
\endproof

\begin{cor} \label{Z2} Let $V \subseteq K^n$ be an $m$-dimensional subspace of $K^n$ and $F$ a multilinear $(n,g)$-form over $K$. Assume that $m+g-1 > n$ and $g > 1$. Then $V$ contains a basis $\bx_1,\dots,\bx_m$ of vectors such that $F(\bx_1) = \dots = F(\bx_m) = 0$ and
\begin{equation}
\label{Z2_bnd}
H(\bx_i) \leq \sqrt{2} m |\Delta_K|^{\frac{m+1}{2d}} \H(V).
\end{equation}
for each $1 \leq i \leq m$.
\end{cor}

\proof
By Theorem~\ref{siegel_sparse}, $V$ contains a basis $\bx_1,\dots,\bx_m$ of $(n-m+1)$-sparse vectors satisfying~\eqref{Z2_bnd}. Since $m+g-1 > n$, these vectors have no more than $g-1 \geq n-m+1$ nonzero coordinates, and hence $F(\bx_i) = 0$ for each $1 \leq i \leq m$.
\endproof

\begin{cor} \label{Z3} Let $V \subseteq K^n$ be an $m$-dimensional subspace of $K^n$ and $F$ a multilinear $(n,g)$-form over $K$. Assume that $m+g-1 > n$ and $g > 1$. Suppose also that $P(x_1,\dots,x_n) \in K[x_1,\dots,x_n]$ is a polynomial such that $D(P,V) < m$. Then there exists a point $\bz \in V \setminus Z_K(P,V)$ such that $F(\bz) = 0$ and
$$H(\bz) \leq \sqrt{2} m |\Delta_K|^{\frac{m+1}{2d}} \H(V).$$
\end{cor}

\proof
Let $\bx_1,\dots,\bx_m$ be the $(n-m+1)$-sparse basis for $V$ guaranteed by Theorem~\ref{siegel_sparse}. Since $m+g-1 > n$, we see that $g > n-m+1$ and hence $F(\bx_i) = 0$ for each $1 \leq i \leq m$, by the same argument as in the proof of Corollary~\ref{Z2}. Since $D(P,V) < m$, at least one of these vectors is not in $Z_K(P,V)$. Call this vector $\bz$, and the result follows.
\endproof
\smallskip

\proof[Proof of Theorem~\ref{mf_zeros}]
The theorem now follows by combining Corollaries~\ref{Z2} and~\ref{Z3}. As such, it is a consequence of the ``sparse" Siegel's lemma (Theorem~\ref{siegel_sparse}).
\endproof
\smallskip

The case $g=1$ of linear forms has to be considered separately: this is just a simple case of Theorem~1.4 of~\cite{lf:null}, which we present here in a simplified form. Suppose that
$$F(x_1,\dots,x_n) = \sum_{i=1}^n a_i x_i \in K[x_1,\dots,x_n],$$
and $V \subseteq K^n$ is an $m$-dimensional subspace.

\begin{lem} \label{linear} Let $V \subseteq K^n$ be an $m$-dimensional subspace with $2 \leq m \leq n$. Then there exists $\bo \neq \bz \in V$ such that $F(\bz) = 0$ and 
\begin{equation}
\label{lin_bnd-1}
H(\bz) \leq \left( \left( \frac{2}{\pi} \right)^{2r_2} |\Delta_K| \right)^{\frac{1}{2d}} \left( \sqrt{n} \H(V) H(F) \right)^{\frac{1}{m-1}}.
\end{equation}
Further, if $P(x_1,\dots,x_n) \in K[x_1,\dots,x_n]$ is a polynomial such that $D(P,V) < m-1$, then there exists a point $\bz \in V \setminus Z_K(P,V)$ such that $F(\bz) = 0$ and
\begin{equation}
\label{lin_bnd-2}
H(\bz) \leq \sqrt{n} \left( \left( \frac{2}{\pi} \right)^{2r_2} |\Delta_K| \right)^{\frac{m-1}{2d}} \H(V) H(F).
\end{equation}
\end{lem}

\proof
Define
$$U(F) = \left\{ \bx \in K^n : F(\bx)=0 \right\},$$
then, by~\eqref{duality_p},
\begin{equation}
\label{ht_UF}
\H(U(F)) = \H(F) \leq \sqrt{n} H(F),
\end{equation}
where $\H(F)$ is the $\H$ height applied to the coefficient vector of $F$, which has $n$ coordinates. This is an $(n-1)$-dimensional subspace of $K^n$ and $\dim V \geq 2$, hence 
$$\ell := \dim (V \cap U(F)) \geq m-1 \geq 1.$$
By Theorem~\ref{siegel_lem}, there is a basis $\bx_1,\dots,\bx_{\ell} \in V \cap U(F)$ such that
\begin{eqnarray}
\label{ht_UF-1}
\prod_{i=1}^{\ell} H(\bx_i) & \leq & \left( \left( \frac{2}{\pi} \right)^{2r_2} |\Delta_K| \right)^{\frac{\ell}{2d}} \H(V \cap U(F)) \nonumber \\
& \leq & \left( \left( \frac{2}{\pi} \right)^{2r_2} |\Delta_K| \right)^{\frac{\ell}{2d}} \H(V) \H(U(F)),
\end{eqnarray}
by Lemma~\ref{intersect}. Then~\eqref{lin_bnd-1} follows by combining~\eqref{ht_UF} with~\eqref{ht_UF-1} and taking $\bz$ to be the vector with smallest height among $\bx_1,\dots,\bx_{\ell}$. Since $D(P,V) < m-1$, at least one of these vectors is not in $Z_K(P,V)$, and this implies~\eqref{lin_bnd-2}.
\endproof

\bigskip
\noindent
{\bf Acknowledgement:} We wish to thank Pavel Guerzhoy for some helpful remarks. We are also grateful to the anonymous referee for many thoughtful suggestions that improved the quality of presentation.

\bibliographystyle{plain}  

\end{document}